\documentclass[11pt]{article}

\usepackage{amssymb}
\usepackage{amsmath,amscd,amssymb,amsfonts,graphicx,color}

\newtheorem{theorem}{Theorem}
\newtheorem{definition}{Definition}
\newtheorem{corollary}{Corollary}
\newtheorem{remark}{Remark}
\newtheorem{lemma}{Lemma}

\newtheorem{proposition}{Proposition}
\def\proof{{\bf Proof.} }
\def\endproof{\hfill $\Box$}

\usepackage{amsmath,amscd,amssymb,amsfonts,graphicx,color}
\usepackage{verbatim}
\usepackage[affil-it]{authblk}

\def\R     {\ensuremath{\mathbb R}}
\def\N     {\ensuremath{\mathbb N}}

\begin{document}

\def\C{\mathbb{C}}
\def\H{\mathcal{H}}
\def\R{\mathbb{R}}
\def\N{\mathbb{N}}
\def\Z{\mathbb{Z}}
\def\a{\alpha}
\def\b{\beta}
\def\g{\gamma}




\date{\today}
\title{A unified matrix approach to the representation of Appell
polynomials\thanks{{This work was supported in part by {\it FEDER}
funds through {\it COMPETE}--Operational Programme Factors of Competitiveness
(``Programa Operacional Factores de Competitividade'') and by Portuguese funds
through the {\it Center for Research and Development in Mathematics and
Applications}  and the Portuguese Foundation for Science
and Technology  (``FCT--Funda\c{c}\~{a}o para a Ci\^{e}ncia e a Tecnologia''),
within project  PEst-C/MAT/UI4106/2011 with COMPETE number
FCOMP-01-0124-FEDER-022690.}}}

\author[a]{Lidia Aceto}   
\affil[a]{Department of Mathematics,
University of Pisa, Italy}

\author[b]{Helmuth R. Malonek}
\affil[b]{Department of Mathematics,
University of Aveiro, Portugal}

\author[c]{Gra\c{c}a Tomaz}
\affil[c]{Department of Mathematics, Polytechnic Institute of Guarda, Portugal}
 
\maketitle

\begin{abstract}
In this paper we propose a unified approach to matrix representations
of different types of Appell polynomials. This approach is based on the creation matrix -
a special matrix which has only the natural numbers as entries and is closely
related to the well known Pascal matrix. By this means we stress the
arithmetical origins of Appell polynomials. The approach also allows to derive,
in a simplified way, the properties of Appell polynomials by using only matrix
operations.
\end{abstract}

\smallskip




\section{Introduction}
In the last years, the interest in Appell polynomials
and their applications in different fields has significantly
increased. As recent applications of Appell polynomials in fields
like probability theory and statistics we mention \cite{Ans} and
\cite{Sal}. Generalized Appell polynomials as tools for approximating
3D-mappings were introduced for the first time in \cite{Mal2} in combination
with Clifford Analysis methods.

As another new example we mention representation theoretic results like those of
\cite{Bra} and \cite{Lav}, that gave evidence to the central role of
Appell polynomials as orthogonal polynomials. Representation theory is also
the tool for their application in quantum physics as explained in \cite{wein}.

Starting from Appell polynomials, but in the general framework of
noncommutative Clifford algebras, one can find more traditionally
motivated operational approaches to  generalize Laguerre, Gould-Hopper and
Chebyshev polynomials in the recent papers \cite{Cac1, Cac3, Cac2}.

At the same time other authors were concerned with
finding new characterizations of Appell polynomials themselves 
through new approaches. To quote
some of them we mention, for instance, the novel approach developed
in \cite{YangY} which makes use of the generalized Pascal functional
matrices introduced in \cite{YangM} and also a new characterization
based on a determinantal definition proposed in \cite{Cos}. Both of
them have permitted to derive some properties of Appell polynomials
by employing only linear algebra
tools and to generalize some classical Appell polynomials.\\

In this paper we concentrate on an unifying tool for representing Appell
polynomials in one real variable in matrix form. The matrix structure of
Appell polynomials is based on one and the
same simple constant matrix $H$ defined by
\begin{equation} \label{H}
(H)_{ij}=\left\{
    \begin{array}{ll}
      i, & \quad  i=j+1 \\
      0, & \quad  \hbox{otherwise,}  \qquad i,j=0,1,\ldots,m.
    \end{array}
  \right.
\end{equation}
Needless to emphasize the simplicity of the structure of $H$: it is a subdiagonal
matrix which contains as nonzero entries only the sequence of natural numbers. In the context
of this paper $H$ is called \textit{creation matrix} because it is the matrix from which
all the special types of Appell polynomials will be created. Our approach has the advantage of reducing the class of Appell polynomials to its
most fundamental arithmetical origins.\\

 At the same time the present paper can be considered as a far reaching
generalization of the papers \cite{Ace1,Ace2, Ace3} where the authors
concentrated their attention on the Pascal matrix and highlighted its
connections with matrix representations of other special polynomials like,
for instance, Bernoulli and Bernstein polynomials. The common Appell polynomial
background was not a subject of their concern.\\

The paper is organized as follows. In Section~$\ref{sec:1}$ the
concept of Appell polynomials and the generating functions of some
particular cases are recalled. Section~$\ref{sec:2}$ explains
how to obtain a matrix representation of the referred
polynomials. In Section~$\ref{sec:3}$ we prove some properties of
Appell polynomials by employing only matrix operations. We also
observe that the approach works for those polynomials $p_n(x)$ which
are not Appell in the strong sense, i.e., where the degree is not
equal to $n,$ using as an example the Genocchi polynomials. Finally,
in Section~$\ref{sec:5}$ some conclusions are presented.

\section{Classical equivalent characterizations of Appell polynomials}
\label{sec:1}
In 1880 Appell introduced in \cite{App} sequences of polynomials
$\{p_n(x)\}_{n\geq 0}$ satisfying the relation
\begin{eqnarray} \label{Apdiff1}
\frac{d}{dx}p_n(x) = n \, p_{n-1}(x), \qquad n=1,2,\dots,
\end{eqnarray}
in which
\begin{eqnarray} \label{Apdiff2}
p_0(x)=c_0, \qquad \quad c_0 \in \R \setminus \{0 \}.
\end{eqnarray}
\begin{remark}
Some authors use another definition for $\{p_n(x)\}_{n\geq
0},$ where the factor $n$ on the right hand side of
(\ref{Apdiff1}) is omitted. However, considering that the prototype of such
sequences is the monomial power basis, i.e., $\{x^n\}_{n\geq 0},$ we
prefer to consider Appell's original definition.
\end{remark}
From (\ref{Apdiff1})-(\ref{Apdiff2}) it can easily be checked that
the polynomials which form an Appell set are $n$-degree polynomials
having the following form:
\begin{eqnarray*}  \label{explicAp}
p_0(x) &=& c_0  \\
p_1(x)&=& c_1+c_0 \, x   \nonumber \\
p_2(x)&=& c_2+2 \,c_1 \,x +c_0 \, x^2  \nonumber  \\
p_3(x)&=& c_3+3\, c_2 \, x+3 \,c_1 \,x^2+c_0 \, x^3 \nonumber  \\
&\vdots& \nonumber
\end{eqnarray*}
This can be written in a compact form as follows:
\begin{eqnarray}\label{compacAp}
p_n(x)= \sum_{k=0}^n \left(
                       \begin{array}{c}
                         n \\
                         k \\
                       \end{array}
                     \right) c_{n-k} \, x^k,  \quad n=0,1,\dots, \qquad  c_0
\neq 0.
\end{eqnarray}
In addition,  if $f(t)$ is any convergent power series on the whole real line
with Taylor expansion given by
\begin{eqnarray} \label{ef}
f(t)=\sum_{n=0}^{+\infty} c_n \, \frac{t^n}{n!},  \qquad  f(0)\neq
0,
\end{eqnarray}
then Appell sequences can also be defined by means of the corresponding
{\em generating function} (see  \cite{Boa}),
\begin{eqnarray} \label{genfun}
G(x,t) \equiv f(t) \, e^{xt} = \sum_{n=0}^{+\infty} p_n(x) \, \frac{t^n}{n!}.
\end{eqnarray}
By appropriately choosing the function $f(t),$ many of the classical
polynomials can be derived. In particular, we get \footnote{The
expression of functions $f(t)$ for Bernoulli and Euler polynomials
can be found in \cite{Rai}, the one of the monic Hermite polynomials
in \cite{Boa}  and that of Laguerre polynomials in \cite{Erd}.}
\begin{itemize}
\item the \emph{monomials} $\{x^n\}_{n\geq 0}$ when
\begin{eqnarray*}
f(t)=1;
\end{eqnarray*}
\item the \emph{Bernoulli polynomials} $\{B_n(x)\}_{n \ge 0}$ when
\begin{eqnarray*} \label{genBER}
f(t) &=&\frac{t}{e^t-1} =\left(\frac{\sum_{n=0}^{+\infty}
\frac{t^n}{n!}-1}{t}\right)^{-1} =\left(E_{1,2}(t)\right)^{-1},
\end{eqnarray*}
where
\begin{eqnarray}\label{mlf}
 E_{\a,\b}(z)=\sum_{n=0}^{+\infty}\frac{z^n}{\Gamma(\a
n+\b)}, \qquad \a,
\b >0,
\end{eqnarray}
is the two parameter function of the Mittag-Leffler type and $\Gamma(z)$ the
Gamma function \cite{Pod};
\item the \emph{Euler polynomials} $\{E_n(x)\}_{n \ge 0}$ when
\begin{eqnarray*} \label{genEUL}
f(t) =\frac{2}{e^t+1};
\end{eqnarray*}
\item the \emph{monic Hermite polynomials} $\{\widehat{H}_n(x)\}_{n \ge 0}\equiv
\{2^{-n}H_n(x)\}_{n\geq 0}$, with $H_n(x)$
 the classical Hermite polynomials, when
\begin{eqnarray*} \label{genHER}
f(t)&=&e^{-\frac{t^2}{4}} =E_{1,1}\left(-\frac{t^2}{^4}\right),
\end{eqnarray*}
(see (\ref{mlf}));
\item the \emph{``modified'' generalized Laguerre polynomials} $\{(-1)^n \, n!
\, L_n^{(\alpha-n)}(x)\}_{n \ge 0}$,
$\alpha>-1,$ when
\begin{eqnarray*} \label{genLAG}
f(t) = (1-t)^\alpha.
\end{eqnarray*}
\end{itemize}
Some others families of special polynomials do not seem to be Appell
sets since they are usually defined by a different type of
generating function. However, Carlson in \cite{Car} highlights that
sometimes it is possible to transform a given polynomial sequence into
one of Appell type by suitable changes of variables.  This is
the case, for example, of the \emph{Legendre polynomials}
$\{P_n(x)\}_{n \ge 0}$ defined on the interval $(-1,1)$ whose
generating function is \footnote{The generating functions of $P_n,
T_n, U_n$ here reported can be found in \cite{Rai}.}
\begin{eqnarray*} \label{genLEG0}
J_0 \left(t \, \sqrt{1-x^2}\right)  \, e^{xt} = \sum_{n=0}^{+\infty}
P_n(x) \, \frac{t^n}{n!},
\end{eqnarray*}
where
\begin{eqnarray*} \label{Bessel}
J_0(y)= \sum_{n=0}^{+\infty} (-1)^n \, \frac{y^{2n}}{2^{2n} \,
(n!)^2}
\end{eqnarray*}
denotes the Bessel function of the first kind and index $0.$
However, setting 
\begin{eqnarray} \label{transf0}
x=\frac{z}{\sqrt{z^2+1}}, \qquad t=\tau \sqrt{z^2+1},
\end{eqnarray}
we get
\begin{eqnarray*} \label{genLEG}
J_0 \left(\tau \right)  \, e^{z \tau} = \sum_{n=0}^{+\infty}
P_n\left(\frac{z}{ \sqrt{z^2+1}}\right) \, \sqrt{(z^2+1)^n}
\,\frac{\tau^n}{n!},
\end{eqnarray*}
which is in the form (\ref{genfun}). \\
Similarly, it can be verified that also the Chebyshev polynomials, both
of the first and of the second kind, are Appell sequences. In fact, the
\emph{Chebyshev polynomials of the first kind} $\{T_n(x)\}_{n \ge 0}$ have as
generating function
\begin{eqnarray} \label{algenCHE1}
\cosh\left(t \, \sqrt{x^2-1}\right)  \, e^{xt} =
\sum_{n=0}^{+\infty} T_n(x) \, \frac{t^n}{n!}, \qquad x \in (-1,1),
\end{eqnarray}
and the \emph{Chebyshev polynomials of the second kind}
$\{U_n(x)\}_{n \ge 0}$ have as generating function
\begin{eqnarray} \label{algenCHE2}
\frac{\sinh\left(t \, \sqrt{x^2-1}\right)}{\sqrt{x^2-1}} \, e^{xt} =
\sum_{n=0}^{+\infty} U_n(x) \, \frac{t^{n+1}}{(n+1)!}, \qquad x \in (-1,1).
\end{eqnarray}
Therefore, by using the transformations (\ref{transf0}) in (\ref{algenCHE1}) as well as in (\ref{algenCHE2})  we
get
\begin{eqnarray*}
\cosh\left(i \, \tau \right)  \, e^{z \tau} = \sum_{n=0}^{+\infty}
\, T_n\left(\frac{z}{ \sqrt{z^2+1}}\right) \, \sqrt{(z^2+1)^n} \,
\frac{\tau^n}{n!},
\end{eqnarray*}
respectively
\begin{eqnarray*}
\frac{\sinh\left(i \, \tau \right)}{i}  \, e^{z \tau} =
\sum_{n=0}^{+\infty} \, U_n\left(\frac{z}{ \sqrt{z^2+1}}\right) \,
\sqrt{(z^2+1)^n} \, \frac{\tau^{n+1}}{(n+1)!}.
\end{eqnarray*}
Consequently, by taking into account the trigonometric identities
\begin{eqnarray*}
\cosh(i\, \tau) &=& \frac{e^{i\, \tau} + e^{-i\, \tau}}{2} = \cos \tau,  \\
\frac{\sinh(i\,\tau)}{i} &=& \frac{e^{i\, \tau} - e^{-i\, \tau}}{2\,
i} = \sin \tau,
\end{eqnarray*}
it follows that the generating functions of the Chebyshev
polynomials of the first and of the second kind are
\begin{eqnarray*}
\cos \tau  \, e^{z \tau} = \sum_{n=0}^{+\infty} \,
T_n\left(\frac{z}{ \sqrt{z^2+1}}\right) \, \sqrt{(z^2+1)^n}\,
\frac{\tau^n}{n!} \label{genCHE1}
\end{eqnarray*}
and
\begin{eqnarray*}
\hbox{sinc}\; \tau \, e^{z \tau} = \sum_{n=0}^{+\infty} \,
\frac{1}{n+1}\, U_n\left(\frac{z}{ \sqrt{z^2+1}}\right) \,
\sqrt{(z^2+1)^n} \, \frac{\tau^n}{n!}, \label{genCHE2}
\end{eqnarray*}
respectively, where, as usual, $\hbox{sinc}\; \tau=\sin \tau/\tau.$ \\
Summarizing, by using the changes of variables  (\ref{transf0}) we
get
\begin{itemize}
\item  the \emph{``modified'' Legendre polynomials}
$\{\sqrt{(z^2+1)^n}P_n(z/\sqrt{z^2+1})\}_{n\geq 0}$ where
\begin{eqnarray*}
f(\tau) = J_0 \left(\tau \right),
\end{eqnarray*}
is the Bessel function of the first kind and index $0;$
\item the \emph{``modified'' Chebyshev polynomials of the first kind}\\
$\{\sqrt{(z^2+1)^n}T_n(z/\sqrt{z^2+1})\}_{n\geq 0}$ where
\begin{eqnarray*} \label{Che1}
f(\tau)=\cos \tau =E_{2,1}(-\tau^2),
\end{eqnarray*}
(see (\ref{mlf}));
\item the \emph{``modified'' Chebyshev polynomials of the second kind}\\
$\{\frac{1}{n+1} \sqrt{(z^2+1)^n}U_n(z/\sqrt{z^2+1})\}_{n\geq 0}$
when
\begin{eqnarray*} \label{Che2}
f(\tau)=\hbox{sinc}\; \tau =E_{2,2}(-\tau^2).
\end{eqnarray*}
\end{itemize}

\section{Appell polynomials: the matrix approach}\label{sec:2}
As mentioned in the Introduction, our unified matrix approach basically
relies on the properties of the {\it creation matrix} (\ref{H}).
It is worth to observe that it is a nilpotent matrix of degree $m+1,$
i.e.,
\begin{eqnarray} \label{matH}
 H^s= O, \qquad \mbox{for all } s \ge m+1.
\end{eqnarray}
This property is one of the essential ingredients for the unified matrix approach to Appell polynomials that now follows.\\

In order to handle the Appell sequence $\{p_n(x)\}_{n\geq 0}$ in a closed form we introduce
\begin{eqnarray}\label{vetpol} {\bf
p}(x)=[p_0(x)\;\;p_1(x)\;\cdots\;p_m(x)]^T,
\end{eqnarray}
hereafter called {\em Appell vector}.

Due to (\ref{Apdiff1}), the application of the creation matrix (\ref{H}) implies that the Appell vector
satisfies the differential equation
\begin{eqnarray} \label{ode}
\frac{d}{dx}{\bf p}(x) = H \, {\bf p}(x),
\end{eqnarray}
whose general solution is
\begin{eqnarray}   \label{genPasc}
{\bf p}(x) = e^{xH} \, {\bf p}(0)  \equiv P(x)  \, {\bf p}(0),
\end{eqnarray}
with $P(x)$ defined by
\begin{eqnarray}   \label{Px}
(P(x))_{ij} = \left\{
    \begin{array}{cl}
    {i \choose j} \, x^{i-j}, & \quad  \textrm{$i\geq j $} \\
      0, & \quad  \hbox{otherwise,}  \qquad i,j=0,1,\ldots,m.
    \end{array}
  \right.
\end{eqnarray}
The matrix (\ref{Px}) is called {\em generalized Pascal matrix}
because $P(1)\equiv P$ is the lower triangular {\em Pascal matrix}
\cite{Ace, Cal} defined by
\begin{eqnarray}   \label{pascal}
(P)_{ij} = \left\{
    \begin{array}{cl}
    {i \choose j}, & \quad  \textrm{$i\geq j $} \\
      0, & \quad  \hbox{otherwise,}  \qquad i,j=0,1,\ldots,m.
    \end{array}
  \right.
\end{eqnarray}
Notice that $P(0)\equiv I$ is the identity matrix.\\

Consider now the vector of monomial powers
$$ \xi(x)=[1\;\;x\;\cdots\;x^m]^T$$
and the matrix $M$ defined by
\begin{eqnarray} \label{matM}
(M)_{ij}=\left\{\begin{array}{lll}
{i \choose j}\, c_{i-j},& \textrm{$i\geq j $} \\
0, & \textrm{otherwise,}& \quad  \textrm{$i,j=0,1,\ldots,m.$}
\end{array}\right.
\end{eqnarray}
According to (\ref{explicAp}), we have
\begin{eqnarray} \label{Apmat}
 {\bf p}(x)=M\xi(x).
\end{eqnarray}
This relation motivates the following definition.
\begin{definition} \label{transf}
The matrix $M$ defined by (\ref{matM}) is called the \emph{transfer
matrix} for the Appell vector (\ref{vetpol}).
\end{definition}
Obviously,
\begin{eqnarray}\label{matrixm}
{\bf p}(0)=M  {\xi}(0)=  [c_0\;\;c_1\;\cdots\;c_m]^T.
\end{eqnarray}
Therefore, from  (\ref{genPasc}) and (\ref{Apmat}) we conclude that to
obtain the different kinds of Appell polynomials we need to specify
the entries of ${\bf p}(0)$ or, equivalently, of the transfer matrix
$M.$ At this aim, a powerful tool is given by the following result.
 \begin{theorem}    \label{teo1} Let $H$ be the creation matrix defined by (\ref{H}).
If $G(x,t) \equiv f(t)e^{xt}$ is the generating function of an
Appell sequence $\{p_n(x)\}_{n\geq 0}$, then the transfer matrix $M$
is a nonsingular matrix equal to $f(H).$
\end{theorem}
\proof
Since (see (\ref{genfun}))
$$G(x,t)=\sum_{n=0}^{+\infty}p_n(x)\frac{t^n}{n!},$$
setting $x=0,$ we obtain
 $$f(t)=\sum_{n=0}^{+\infty}p_n(0)\frac{t^n}{n!}.$$
Substituting $t$ by $H$ and by taking into account (\ref{ef}) and
(\ref{matH}) we get
$$f(H)=\sum_{n=0}^{m}c_n\frac{H^n}{n!}, \qquad \quad c_0 \neq 0.$$
Furthermore, denoting by ${\bf e}_s, s=0,1,\ldots,m$, the standard
unit basis vector in $\R^{m+1}$ and adopting the convention that
${\bf e}_s={\bf 0}$ whenever $s>m$ (null vector), then the entries
of $f(H)$ are obtained by
\begin{eqnarray} \label{enterH}
(f(H))_{ij}&=&\sum_{n=0}^m\frac{c_n}{n!}{\bf e}_i^TH^n{\bf
e}_j=\sum_{n=0}^m\frac{c_n}{n!}(j+1)^{(n)}{\bf e}_i^T{\bf
e}_{j+n}\\
&=&\sum_{n=0}^mc_n\frac{(j+1)^{(n)}}{n!}\delta_{i,j+n}, \nonumber
\end{eqnarray}
where $\delta_{i,j}$ is the Kronecker symbol and
$(j+1)^{(n)}=(j+1)(j+2)\cdots(j+n)$ is the ascending factorial with
$(j+1)^{(0)}:=1.$
 Thus, $(f(H))_{ij}=0$ if $i<j$, and when $i=j+n$, that is $i\geq j$,
$$(f(H))_{ij}=c_{i-j}\frac{i!}{(i-j)!j!}={i \choose j}c_{i-j}$$
which completes the proof. \endproof

Taking into account the previous theorem, we achieve for the
particular cases referred in Section~$2$ the corresponding
transfer matrices.
\begin{description}
  \item[(i)] For the monomials $\{x^n\}_{0\leq n \leq m}$, $$M=I,$$
 where $I$ is the identity matrix of order $m+1.$
  \item[(ii)] For the Bernoulli polynomials $\{B_n(x)\}_{0\leq n\leq m}$,
\begin{eqnarray} \label{matBern}
 M=\left(\sum_{n=0}^m\frac{H^n}{(n+1)!}\right)^{-1}.
\end{eqnarray}
  \item[(iii)] For the Euler polynomials $\{E_n(x)\}_{0\leq n\leq m}$,
\begin{eqnarray} \label{matEul}
 M=2(e^H+I)^{-1} = 2(P+I)^{-1},
\end{eqnarray}
  where $P$ is the Pascal matrix of order $m+1$ (see
  (\ref{pascal})).
\item[(iv)] For the monic Hermite polynomials
$\{\widehat{H}_n(x)\}_{0\leq n\leq m}$,
\begin{eqnarray} \label{matHer}
M=e^{-\frac{H^2}{4}}=\sum_{n=0}^m\frac{(-1)^nH^{2n}}{2^{2n} n!}.
\end{eqnarray}
In this case, introducing the diagonal matrix
\begin{eqnarray} \label{Dg}
 D[\ell]=\hbox{diag} [\ell^0,\ell^1,\ell^2,\dots,\ell^m],\qquad \ell \neq 0,
\end{eqnarray}
and ${\bf H}(x)=[H_0(x)\;\;H_1(x)\;\cdots\;H_m(x)]^T,$  the vector
of  the classical Hermite polynomials, we get
$$\left(D[2]\right)^{-1}{\bf H}(x)=M\xi(x) \quad \Leftrightarrow
\quad {\bf H}(x)=D[2]M\xi(x).$$
 \item[(v)] For the ``modified'' generalized Laguerre polynomials
 $$\{(-1)^n n!L_n^{(\alpha-n)}(x)\}_{0\leq n\leq m}$$
we have $$ M=\left(I-H\right)^{\alpha}.$$
It is worth noting that, this matrix point of view provides an easy
way to relate $\{(-1)^n n!L_n^{(\alpha-n)}(x)\}_{0\leq n\leq m}$
with the generalized Laguerre polynomials
$\{L_n^{(\alpha)}(x)\}_{0\leq n\leq m}$ \cite{Rai}. In fact,
introducing the vector
$${\bf L}(x)=[L_0^{(\a)}(x)\;\;
L_1^{(\a-1)}(x)\;\cdots\;L_m^{(\a-m)}(x)]^T$$  and the diagonal matrix
$$D_f=\hbox{diag}[0!,1!,2!,\cdots,m!],$$ we obtain (see (\ref{Dg})) $$D[-1]D_f
{\bf L}(x)=(I-H)^{\a}\xi(x)$$
or, equivalently,
\begin{eqnarray}  \label{laggen}
 {\bf L}(x)=\left( D_f \right)^{-1}  D[-1](I-H)^{\a}\xi(x).
\end{eqnarray}
In addition,  the recurrence relation reported in
\cite{Rai}
$$L_n^{(\a)}(x)=L_{n-1}^{(\a)}(x)+L_n^{(\a-1)}(x),\;\;n>0,$$ gives,
successively,
\begin{eqnarray*}
L_1^{(\a-1)}(x)&=&L_1^{(\a)}(x)-L_0^{(\a)}(x)\\
L_2^{(\a-2)}(x)&=&L_2^{(\a-1)}(x)-L_1^{(\a-1)}(x)\\
&=&(L_2^{(\a)}(x)-L_1^{(\a)}(x))-(L_1^{(\a)}(x)-L_0^{(\a)}(x))\\
&=&L_2^{(\a)}(x)-2L_1^{(\a)}+L_0^{(\a)}(x)\\
&\vdots&{}\\
L_m^{(\a-m)}(x)&=&\sum_{n=0}^m(-1)^{m-n}{m \choose
n}L_{n}^{(\a)}(x).
\end{eqnarray*}
Denoting by
$\boldmath{\mathcal{L}}(x)=[L_0^{(\a)}(x)\;\;L_1^{(\a)}(x)\;\cdots\;L_m^{
(\a)} (x)]^T$ the vector of the first $m+1$ generalized Laguerre
polynomials and taking into account (\ref{Px}), we have
$${\bf L}(x)=P(-1) \boldmath{\mathcal{L}}(x).$$
Finally, from (\ref{laggen}) it follows that
$$\boldmath{\mathcal{L}}(x)=P(D_f)^{-1}D[-1](I-H)^\a
\xi(x).$$ In particular, when $\a=0$ we get the ordinary {\em Laguerre
polynomials}.
\item[(vi)] For the modified Legendre polynomials \\
$\{ \sqrt{(z^2+1)^n} P_n(z/\sqrt{z^2+1})\}_{0\leq n\leq m}$ we get
\begin{eqnarray*}\label{legendre}
  M=J_0(H) = \sum_{n=0}^m (-1)^n\frac{H^{2n}}{2^{2n}(n!)^2}.
 \end{eqnarray*}
\item[(vii)] For the modified Chebyshev polynomials of the first
kind \\$\{ \sqrt{(z^2+1)^n} T_n(z/\sqrt{z^2+1})\}_{0\leq n\leq m}$,
 \begin{eqnarray*}\label{chebyshev1}
  M=\cos H =\sum_{n=0}^m (-1)^n\frac{H^{2n}}{(2n)!}.
  \end{eqnarray*}
  \item[(viii)] For the modified Chebyshev polynomials of the second
kind \\
$\{\frac{1}{n+1}\{ \sqrt{(z^2+1)^n}
U_n(z/\sqrt{z^2+1})\}_{0\leq n\leq m}$,
\begin{eqnarray}\label{chebyshev2}
  M \equiv M_{\bf U} = \sum_{n=0}^m (-1)^n\frac{H^{2n}}{(2n+1)!}.
\end{eqnarray}
\end{description}
Of course, considering
$$
z=\frac{x}{\sqrt{1-x^2}}, \quad  \quad   x\in (-1,1),
$$
we obtain the first $m+1$ elements of the classical sequences
$\{P_n(x)\}_{n\geq 0},$ $\{T_n(x)\}_{n\geq 0}$ and
$\{U_n(x)\}_{n\geq 0}.$ Collecting these elements into the vectors
${\bf P}(x),$ ${\bf T}(x),$ and ${\bf U}(x),$ respectively, we get
\begin{eqnarray*}
{\bf P}(x) &=& D[ \sqrt{1-x^2} ]\, J_0(H) D^{-1}[ \sqrt{1-x^2} ]\,
\xi(x), \\
{\bf T}(x)&=&D[ \sqrt{1-x^2} ]\,  \cos H D^{-1}[ \sqrt{1-x^2} ]\, \xi(x), \\
{\bf U}(x)&=&D_{m+1} D[ \sqrt{1-x^2} ]\, M_{\bf U}  D^{-1}[
\sqrt{1-x^2} ]\, \xi(x),
\end{eqnarray*}
where $D_{m+1}=\hbox{diag} [1, 2, \dots, m+1]$.
\begin{remark}
The matrix $M_{\bf U}$ given in (\ref{chebyshev2}) 
 satisfies $H M_{\bf U} = \sin H.$
\end{remark}

\section{Some properties of Appell polynomials} \label{sec:3}
In order to establish several identities of Appell polynomials in a friendly
and unified way, we now use the transfer matrix and some properties
of the generalized Pascal matrix. It is worth mentioning that some of them
were derived in \cite{Cos,Cos1} by using the determinantal approach.

\begin{lemma} \label{Pcsi}
Let $P(x)$ be the generalized Pascal matrix and  $\xi(x)$ the vector containing
the ordinary monomials as previously defined. Then,
$$
\xi(x+y) = P(x) \, \xi(y), \qquad \forall \, x,y \in \R.
$$
\end{lemma}
\proof
The result is a consequence of the binomial theorem. In fact, (see (\ref{Px}))
$$
\left(\xi(x+y)\right)_i= (x+y)^i = \sum_{k=0}^{i} {i \choose k} x^{i-k} \, y^k
=\left(P(x) \, \xi(y)\right)_i.
$$
\endproof \\
According to Carlson \cite[Theorem~1, p. 545]{Car}, it is known that if
$\{p_n(x)\}_{n\geq 0}$ is an Appell sequence, then it satisfies a
binomial theorem of the form
\begin{equation}\label{binomial}
p_n(x+y)=\sum_{k=0}^\infty {n \choose k}p_k(x)y^{n-k},\qquad
n=0,1,\ldots.
\end {equation}
This property is necessary but not sufficient to define an Appell
sequence in the sense of (\ref{Apdiff1})-(\ref{Apdiff2}). In fact,
$\{p_n(x)\}_{n\geq 0}$ verifies (\ref{binomial}) if and only if
\begin{equation}\label{condappell}
\frac{d}{dx}p_n(x)  =np_{n-1}(x),\;\;\;n=0,1,\ldots,
\end{equation}
where the right hand-side is taken to be zero in the case of $n=0.$
This means that the binomial theorem property does not require the
condition that $p_n(x)$ should be exactly of degree $n$ as we
mentioned in the beginning of Section~$\ref{sec:1}$ as
consequence of (\ref{Apdiff1})-(\ref{Apdiff2}).

\begin{remark} The sequence of the {\em Genocchi polynomials}, whose
generating function is \footnote{The expression of the function
$f(t)$ for the Genocchi polynomials can be found in \cite{Liu}.}
\begin{equation}\label{genocchi}
f(t)e^{xt}=\sum_{n=0}^\infty G_n(x)\frac{t^n}{n!}, \qquad
f(t)=\frac{2t}{e^t+1},
\end{equation}
is an example of a sequence that satisfies (\ref{binomial}).
Nevertheless, by virtue of (\ref{condappell}) and (\ref{genocchi}),
it is possible to get the corresponding transfer matrix $M=M_{\bf
G}$ in the same way as in the Section~$\ref{sec:2}$ for the case of
Appell polynomials, but in this case $M$ is singular. Actually, from
Theorem~\ref{teo1},
\begin{eqnarray*}
M_{\bf G}= 2H(e^H+I)^{-1}=2H(P+I)^{-1}.
\end{eqnarray*}
\end{remark}

Now, from  (\ref{binomial}) we can derive the following result.
\begin{theorem} \label{teo2} Let $\{p_n(x)\}_{n\ge 0}$ be an Appell sequence
and $P(x)$ the generalized Pascal matrix defined by (\ref{Px}).
For the corresponding Appell vector we have
\begin{eqnarray} \label{pxy}
{\bf p}(x+y) = P(y) \,  {\bf p}(x), \qquad \forall \, x,y \in \R.
\end{eqnarray}
\end{theorem}
\proof From (\ref{Apmat}) and Lemma~\ref{Pcsi} one has
$${\bf p}(x+y) = M \xi(x+y) = M P(y) \,  \xi(x).$$
The proof is completed by observing that  $P(y)$ and $M$ commute
since they are both functions of the creation matrix $H.$
\endproof

\begin{corollary} Let $\{p_n(x)\}_{n\ge 0}$ be an Appell sequence. Then,
for any constant $a$ and for all $x \in \R,$
the Appell vector of the given sequence satisfies the following relations:
\begin{itemize}
\item[(i)] forward difference:
$$
 {\bf p}(x+1) -  {\bf p}(x) = (P-I) \,  {\bf p}(x);
$$
\item[(ii)] multiplication theorem:
\begin{eqnarray}
{\bf p}(a x) &=& P((a-1)\,x) \,  {\bf p}(x), \label{pa1}\\
{\bf p}(a x) &=& M D[a] \,  {\bf \xi} (x),  \label{pa2}
\end{eqnarray}
where $D[a]$ is defined by (\ref{Dg}).
\end{itemize}
\end{corollary}
\proof
\begin{itemize}
\item[(i)] The result follows from (\ref{pxy}) with $y=1$ and by recalling
that $P(1) \equiv P.$
\item[(ii)] The relation (\ref{pa1}) can be immediately
deduced from (\ref{pxy}) with $y=(a-1)x.$ Concerning (\ref{pa2}), it
follows from (\ref{Apmat}) and by observing that $\xi(ax)=D[a] \,
\xi(x).$
\endproof
\end{itemize}
It is worth noting that (\ref{pa1}) generalizes, for all kinds of
Appell polynomials, the well-known properties for the Bernoulli and
the Euler polynomials (see, e.g., \cite{Abra, Cos, Raabe})
 \begin{eqnarray*}
B_n(ax) &=&   \sum_{i=0}^{n} {n \choose i}
B_i(x) (a-1)^{n-i} x^{n-i},  \\
E_n(ax) &=&
 \sum_{i=0}^{n} {n \choose i}
E_i(x) (a-1)^{n-i} x^{n-i}.
\end{eqnarray*}
In addition, there are some identities involving Appell polynomials
which at a first glance are not equivalent. To provide an example we
consider the following ones involving the Bernoulli polynomials:
$$
B_n(1-x) = (-1)^n B_n(x), \qquad  B_n(1) = (-1)^n B_n(0).
$$
It is trivial to check that the relation on the left hand-side
implies the one on the right, while it is not clear that the
opposite implication also holds true. Similar arguments can be used
about the following relations involving Euler polynomials:
$$
E_n(1-x) = (-1)^n E_n(x), \qquad  E_n(1) = (-1)^n E_n(0).
$$
However, the referred equivalences will become evident from the next
theorem.
\begin{theorem} \label{symmetry}
Let $\{p_n(x)\}_{n\ge 0}$ be an Appell sequence. For
the corresponding Appell vector the following equivalence holds
\begin{eqnarray}
\left( \, {\bf p}(h-x) = D[-1]\, {\bf p}(x), \, \forall \, h,x \in
\R \, \right) \, \Leftrightarrow \, \left(\, {\bf p}(h) = D[-1]\,
{\bf p}(0), \, \forall  \, h \in \R\,\right) \!,
\end{eqnarray}
where $D[-1]$ is defined by (\ref{Dg}).
\end{theorem}
\proof   ($\Rightarrow$) This implication is trivial from the hypothesis with $x=0.$ \\
($\Leftarrow$) Using (\ref{pxy}) and by observing that $P(-x) = D[-1] P(x) D[-1]
$ (see (\ref{Px})) we get
\begin{eqnarray*}
 {\bf p}(h-x) &=& P(-x)  {\bf p}(h)  = D[-1] P(x) D[-1] D[-1]\, {\bf p}(0) \\
&=& D[-1] P(x)  \, {\bf p}(0) =   D[-1]  \, {\bf p}(x).
\end{eqnarray*}
\endproof

It is interesting to ask for the consequences of Theorem~\ref{symmetry}  for 
the coefficients of the Appell polynomials.

\begin{corollary}  Let $\{p_n(x)\}_{n\ge 0}$ be an Appell sequence. For
the Appell vector we get
\begin{eqnarray*}
\left( \, {\bf p}(-x) = D[-1]\, {\bf p}(x),  \quad \forall\,  x
\in \R \, \right) \quad  \Leftrightarrow
\quad \left(\,  c_{2n+1} = 0, \quad n=0,1,\dots\,  \right).
\end{eqnarray*}
\end{corollary}
\proof   From Theorem~\ref{symmetry}, by fixing $h=0$ one has that
$$
\left( \, {\bf p}(-x) = D[-1]\, {\bf p}(x),  \quad \forall\,  x
\in \R \, \right) \quad  \Leftrightarrow
\quad (\,{\bf p}(0) = D[-1]\, {\bf p}(0)\,).
$$
The relation on the right of this equivalence together with the fact
that ${\bf p}(0) =   [c_0\;\;c_1\;\cdots\;c_m]^T$ (see
(\ref{matrixm})) lead to the desired result.
\endproof \\

We observe that the polynomials of Hermite, Legendre and Chebyshev
of the first and of the second kind belong to the Appell class and verify the
equivalence of the previous corollary. This is due to the fact that their
transfer matrices are defined as expansions of even powers of $H$ whose entries
satisfy (see (\ref{enterH}))
$$
{\bf e}^T_i H^{2n} {\bf e}_j = (j+1)^{(2n)} {\bf e}^T_i {\bf
e}_{j+2n} = 0, \qquad i-j \neq 2n.
$$

Among the properties of Appell polynomials we now refer to the ones
stated in Theorem~$11$ of \cite{Cos}, which were proved applying a determinantal approach. Here we
propose an alternative proof based on our matrix approach and making use
of the transfer matrix.
\begin{theorem} Let $\{v_n(x)\}_{n\ge 0}$ and $\{u_n(x)\}_{n\ge 0}$ be two
sequences of Appell polynomials and ${\bf v}(x)$ and ${\bf u}(x)$
their corresponding Appell vectors. Then,
\begin{itemize}
 \item[$(i)$] for all $\lambda, \mu \in \R,$ $\lambda {\bf v}(x)  + \mu {\bf u}(x)$ is an
Appell vector for the sequence  $\{\lambda v_n(x) + \mu
u_n(x)\}_{n\ge 0};$
 \item[$(ii)$] replacing in $v_n(x)$ the powers $x^0,$
$x^1,\dots,x^n$ by $u_0(x),$ $u_1(x),\dots,u_n(x)$ and denoting
the resulting polynomial by $w_n(x)$, the vector ${\bf
w}(x)=[w_0(x) \, w_1(x)\, \dots \, w_m(x)]^T$ is the Appell vector
of $\{w_n(x)\}_{n\ge 0}.$
\end{itemize}
\end{theorem}
\proof  Let  $M_v$ and $M_u$ be the transfer matrices for ${\bf
v}(x)$ and ${\bf u}(x),$ respectively, i.e.,
$$
{\bf v}(x) = M_v \xi(x), \qquad {\bf u}(x) = M_u \xi(x).
$$
Then,
\begin{itemize}
 \item[$(i)$] $\lambda {\bf v}(x)  + \mu {\bf u}(x) = ( \lambda M_v + \mu M_u)
\, \xi(x), $
 \item[$(ii)$] ${\bf w}(x) = M_v {\bf u}(x).$
\end{itemize}
To prove that $\lambda {\bf v}(x) + \mu {\bf u}(x)$ and ${\bf w}(x)$
are Appell vectors, we need to check if they satisfy the relation
(\ref{ode}). The result follows by recalling that $M_v$ and $M_u$
are both functions of $H$ (see Theorem~\ref{teo1}). Consequently,
they commute with $H$ as well any of their linear combination.
\endproof \\

In particular, the last theorem allows us to obtain in a straightforward way
some classes of Appell polynomials recently
introduced in \cite{Kha}. In fact, the transfer matrix
\begin{itemize}
 \item[$(i)$] for the $2$-iterated Bernoulli  polynomials $\{B^{[2]}_n
(x)\}_{0\le n\le m}$ is (see (\ref{matBern}))
$$M= \left(\sum_{n=0}^m\frac{H^n}{(n+1)!}\right)^{-2};$$
 \item[$(ii)$] for the $2$-iterated Euler polynomials $\{E^{[2]}_n(x)\}_{0\le n\le m}$ is (see (\ref{matEul}))
$$M=4(e^H+I)^{-2}\equiv 4(P+I)^{-2};$$
 \item[$(iii)$] for the Bernoulli-Euler polynomials $\{\! ~_EB_n (x)\}_{0\le n\le m}$ (or
Euler-Bernoulli polynomials $\{\! ~_BE_n (x)\}_{0\le n\le m}$)  is
(see  (\ref{matBern}) and (\ref{matEul}))
$$M= 2 \left( (P+I) \, \sum_{n=0}^m\frac{H^n}{(n+1)!}\right)^{-1}.$$
\end{itemize}
Some other properties of Appell polynomials can be obtained by making
use of the inverse of the transfer matrix. To achieve them, we
recall that (see Theorem~\ref{teo1} and (\ref{ef}))
$$ M = \sum_{k=0}^m c_k \frac{H^k}{k!}, \qquad \quad c_0 \neq 0.$$
Setting
\begin{equation} \label{matinv}
M^{-1} = \sum_{k=0}^m \gamma_k \frac{H^k}{k!},
\end{equation}
we have
\begin{equation} \label{coefN}
 \gamma_0 = \frac{1}{c_0}, \qquad \gamma_k =- \frac{1}{c_0} \sum_{s=0}^{k-1}
{k \choose s}\, c_{k-s} \gamma_s, \qquad k=1,2,\dots,m.
\end{equation}
In fact,
\begin{eqnarray*}
I&=&M
M^{-1}=\left(\sum_{k=0}^{m}c_k\frac{H^k}{k!}\right)\left(\sum_{r=0}^
{m} \gamma_r \frac{H^r}{r!}\right)
=\sum_{n=0}^{m}\left(\sum_{k+r=n}c_k \gamma_r\frac{H^{k+r}}{k!r!}\right)\\
 &=&\sum_{n=0}^{m}\left(\sum_{r=0}^n
\frac{n!c_{n-r} \gamma_r}{(n-r)!r!}\right)\frac{H^n}{n!}
=\sum_{n=0}^{m} \left(\sum_{r=0}^n {n \choose r}c_{n-r}\gamma_r
\right)\frac{H^n}{n!}.
\end{eqnarray*}
Consequently, (\ref{Apmat}) implies that
$$ M^{-1} {\bf p}(x) = \xi(x)$$
or, equivalently,
$$
\sum_{k=0}^n {n \choose k}\, \gamma_{n-k} \, p_k(x) = x^n, \qquad
n=0,1,\dots,m,
$$
from which we deduce a \textit{general recurrence relation} for Appell
polynomials:
\begin{eqnarray*}
p_n(x) = \frac{1}{\gamma_0} \left(x^n - \sum_{k=0}^{n-1}  {n \choose k}
\gamma_{n-k} \, p_k(x) \right), \qquad n=0,1,\dots.
\end{eqnarray*}

By taking into account (\ref{matinv}), it is an easy matter to notice from (\ref{matBern}),
(\ref{matEul}) and (\ref{matHer}), that
\begin{itemize}
\item for the Bernoulli polynomials:
$$
\gamma_k = \frac{1}{k+1}, \qquad k=0,1,\dots, m;
$$
\item for the Euler polynomials:
$$
\gamma_0=1, \quad \gamma_k = \frac{1}{2}, \qquad k=1,\dots, m;
$$
\item for the monic Hermite polynomials:
$$
\gamma_k = \left\{\begin{array}{ll}
             \displaystyle \frac{1}{2^k}, & k \mbox{ even }  \\
             0, & k \mbox{ odd} \\
           \end{array}\right., \qquad k=0,1,\dots,m.
$$
\end{itemize}
Furthermore, for the {\em generalized Euler polynomials} introduced
in \cite{Cos} we have
$$
\gamma_0=1, \quad \gamma_k = {\bar \gamma}, \qquad k=1,\dots, m,
$$
which leads to
\begin{eqnarray*}
M^{-1} = I + \sum_{k=1}^m {\bar \gamma}\, \frac{H^k}{k!} =  (1-{\bar
\gamma}) I + \sum_{k=0}^m {\bar \gamma}\, \frac{H^k}{k!} = (1-{\bar
\gamma}) \,I +  {\bar \gamma}\,  P.
\end{eqnarray*}
Knowing the relationship between the coefficients of $M$ and its
inverse, we can prove the following result which relates the coefficients of an
Appell polynomial with those of the general recurrence relation.
\begin{proposition}
The elements of the sets $\{c_n\}_{0 \le n \le m}$ and
$\{\gamma_n\}_{0\le n \le m}$ characterizing the transfer matrix $M$
and its inverse, respectively, satisfy the following equivalences:
 $$  c_{2j+1} = 0 \quad  \Leftrightarrow \quad \gamma_{2j+1} = 0, \qquad \quad j=0,1,\dots,\frac{m-1}{2}.$$
\end{proposition}
\proof We proceed by induction on $j.$ If $j=0,$ the statement is
verified directly using (\ref{coefN}). Next, let us suppose that it
is true for $j-1$ and we prove that it is true also for $j.$
From (\ref{coefN})
\begin{eqnarray*}
\gamma_{2j+1} = &-&  \frac{1}{c_0} \sum_{s=0}^{2j} {2j+1 \choose
s}\,
c_{2j+1-s} \gamma_s = - \frac{1}{c_0} c_{2j+1} \gamma_0 +\\
&-&  \frac{1}{c_0}   \sum_{\substack{1\le s\le 2j-1 \\ s \,\, odd}}
{2j+1 \choose s}\, c_{2j+1-s} \gamma_s - \frac{1}{c_0}
 \sum_{\substack{2\le s\le 2j \\ s \,\, even}}\, c_{2j+1-s} {2j+1 \choose
s}\,\gamma_s.
\end{eqnarray*}  
By using the induction hypothesis the last two
sums vanish, and this completes the proof. \endproof

\section{Conclusion} \label{sec:5}
For almost all classical polynomials defined in the ordinary way as,
for instance, in \cite{Rai}, the corresponding generating
functions are well known. In some cases, like Bernoulli polynomials
or Euler polynomials, the usual generating functions are already
given in a form that
reveals their property of being Appell polynomials due to the
inclusion of the exponential function (see (\ref{genfun}) and
\cite{Boa}).

But this is not the case for all classical polynomials, like Legendre or
Chebyshev (both of the first and second kind) polynomials. In these cases,
this paper shows how the Appell polynomial nature can be disclosed by some
substitution in a way that they can be treated as such. In this sense the paper
tries to call attention also to important and well known polynomials that normally are
not known as Appell polynomials.

Being the central ingredients of the presented unified matrix approach to Appell
polynomials, the roles of the creation matrix $H$ as well as of the transfer matrix $M$ are studied.

Furthermore, the paper confirmed the effectiveness of the unified
matrix representation by showing that some new types of recently
introduced Appell polynomials can immediately be deduced. 

Finally, the special role of the transfer matrix is also stressed and
advantageously used for deriving, in an easy and compact way, the relationship between
the coefficients of Appell polynomials and their general recurrence relations.








\end{document}